\documentclass[12pt]{amsart}
\usepackage{amsmath,amssymb,latexsym}

\theoremstyle{plain}
\newtheorem{Thm}{Theorem}
\newtheorem{Cor}{Corollary}

\newtheorem{Def}{Definition}

\input{epsf}

\theoremstyle{Def}

\begin{document}

\title
{On The Kauffman Skein Modules}

\author{Jianyuan K. Zhong}

\address{Department of Mathematics \& Statistics\\
   Louisiana Tech University\\
    Ruston, LA 71272}
\email{zhong@coes.latech.edu}
\author{Bin Lu}

\address{Department of Mathematics \\
   University of Arizona\\
    Tucson, AZ 85721}
\email{binlu@math.arizona.edu}

\date{July 25, 2001}


\begin{abstract}
Let $k$ be a subring of the field of rational functions in
$\alpha, s$ which contains $\alpha^{\pm 1}, s^{\pm 1}$. Let $M$ be a compact oriented $3$-manifold, and let $K(M)$ denote the Kauffman skein module of $M$ over $k$. Then $K(M)$ is the free $k$-module generated
by isotopy classes of framed links in $M$ modulo the Kauffman
skein relations. In the case of  $k={\mathbb Q}(\alpha, s)$, the
field of rational functions in $\alpha, s$, we give a basis for
the Kauffman skein module of the solid torus and a basis for the
relative Kauffman skein module of the solid torus with two points
on the boundary. We then show that $K(S^{1} \times S^2)$ is
generated by the empty link, i.e., $K(S^{1} \times S^2)\cong k$.
\end{abstract}

\maketitle
\section{Introduction}

In a previous paper \cite{ZG99}, Gilmer and the
first author investigate the Homflypt skein module of $S^1\times
S^2$. Here we follow similar approaches in studying the Kauffman
skein modules. The Homflypt skein and Kauffman skein are closely
related. The relative Homflypt skein module of a cylinder with
$2n$ framed points ($n$ input points and $n$ output points) is a
geometric realization of the $n$th Hecke algebra $H_n$ of type A,
and the relative Kauffman skein module of a cylinder with $2n$
points is a geometric realization of the Birman-Murakami-Wenzl
algebra $K_n$. Over ${\mathbb Q}(\alpha, s)$, the field of
rational functions in $\alpha, s$, Beliakava and Blanchet
\cite{BB2000} have given a canonical projection $\pi_n: K_n\to
H_n$ and a multiplicative homomorphism $s_n: H_n\to K_n$. Their
work reveals close connection between $H_n$ and $K_n$. Our work on
Kauffman skein modules is motivated by Gilmer and Zhong's previous
work \cite{ZG99} on Homflypt skein modules and Beliakava and
Blanchet's work on the Birman-Murakami-Wenzl algebra $K_n$
\cite{BB2000}.

Let $k$ be an integral domain containing the invertible elements
$\alpha $ and $s$. We assume that $s-s^{-1}$ is invertible in $k$.

By a {\it framed oriented link} we mean a link equipped with a
string orientation together with a nonzero normal vector field up
to homotopy. By a {\it  framed link} we mean an unoriented framed
link. The links described by figures in this paper will be
assigned the vertical framing which points towards the reader.

Let $M$ be a smooth, compact and oriented $3$-manifold.
\begin{Def}
{\bf The Homflypt skein module} of $M$, denoted by $S(M)$, is the
$k$-module freely generated by isotopy classes of framed oriented
links in $M$ including the empty link modulo the Homflypt skein
relations given in the following figure:
$$
\raisebox{-3mm}{\epsfxsize.3in\epsffile{left.ai}}\quad
        -\quad\raisebox{-3mm}{\epsfxsize.3in\epsffile{right.ai}}\quad
=\quad (\ s - \
s^{-1})\quad\raisebox{-3mm}{\epsfxsize.3in\epsffile{parra.ai}}\quad
,
$$
$$
\raisebox{-3mm}{\epsfxsize.35in\epsffile{framel.ai}}
        \quad=\quad \alpha\quad \raisebox{-3mm}{\epsfysize.3in\epsffile{orline.ai}}\quad  ,
$$
$$
L\ \sqcup \raisebox{-2mm}{\epsfysize.3in\epsffile{unknot.ai}}\quad
                                        =\quad {\dfrac{\alpha-\alpha^{-1}} {\ s - \ s^{-1}}}\quad L\quad.
$$
The last relation follows from the first two if $L$ is nonempty.
\end{Def}
\noindent {\bf Remark:} The definition here  with variables
$\alpha$, $s$ is a version of the original definition with three
variables $x$, $v$ and $s$ in \cite{ZG99} by taking $x=1$,
$v={\alpha}^{-1}$ and using the same $s$. The previous results on
Homflypt skein modules will be carried over under this
specialization.
\begin{Def}
{\bf The relative Homflypt skein module.} \\ Let $X=\{x_1,
x_2,\cdots, x_{n}\}$ be a finite set of framed points oriented
negatively (called input points) in the boundary $\partial M$, and let
$Y=\{y_1, y_2,\cdots, y_{n}\}$ be a finite set of framed points
oriented positively (called output points) in
$\partial M$. Define the relative skein module $S(M, X, Y)$ to be
the $k$-module generated by relative framed oriented links in
$(M,\partial M)$ such that $L \cap \partial M=\partial L=\{x_i,
y_i\}$ with the induced framing and orientation, considered up to
an ambient isotopy fixing $\partial M$, and quotiented by the Homflypt
skein relations.
\end{Def}

In the cylinder $D^2\times I$, let $X_n$ be a set of $n$ distinct
input framed points on a diameter $D^2\times \{1\}$ and $Y_n$ be a
set of $n$ distinct output framed points on a diameter $D^2\times
\{0\}$, it is a well-known result \cite{AM98} that the relative
Homflypt skein module $K(D^2\times I, X_n\amalg Y_n)$ is
isomorphic to the {\it $n$th Hecke algebra} $H_n$, which is the
quotient of the braid group algebra $k[B_n]$ by the Homflypt skein
relations. In section 4,
we will also use the work by Gilmer and Zhong \cite{ZG99} on $S(S^1\times D^2, A, B)$, which is the relative Homflypt skein module of the solid torus with $A$ an input point and $B$ an output point on the boundary of $S^1\times D^2$. 
\begin{Def}
{\bf The Kauffman skein module} of $M$, denoted by $K(M)$, is the
$k$-module freely generated by isotopy classes of framed links in
$M$ including the empty link modulo the Kauffman skein relations
given in the following figure:
$$
\raisebox{-3mm}{\epsfxsize.3in\epsffile{left1.ai}}\quad
        -\quad\raisebox{-3mm}{\epsfxsize.3in\epsffile{right1.ai}}\quad
=\quad (\ s - \
s^{-1})\Bigl(\quad\raisebox{-3mm}{\epsfxsize.3in\epsffile{parra1.ai}}-\quad\raisebox{-3mm}{\epsfxsize.3in\epsffile{parra2.ai}}\quad\Bigr)\quad
,
$$
$$
\raisebox{-3mm}{\epsfxsize.35in\epsffile{framel1.ai}}
        \quad=\quad \alpha \quad \raisebox{-3mm}{\epsfysize.3in\epsffile{orline1.ai}}\quad  ,
$$
$$
L\ \sqcup
\raisebox{-2mm}{\epsfysize.3in\epsffile{unknot1.ai}}\quad =\quad
({\dfrac{\alpha-\alpha^{-1}} {\ s - \ s^{-1}}}+1)\quad L\quad.
$$
The last relation follows from the first two when $L$ is nonempty.
\end{Def}

Similarly, we give the definition of the relative Kauffman skein module as the following.

\begin{Def}
{\bf The relative Kauffman skein module}. \\ Let $X=\{x_1,
x_2,\cdots, x_{n}, y_1, y_2,\cdots, y_{n}\}$ be a finite set of
$2n$ framed points in the boundary $\partial M$. Define the
relative skein module $K(M, X)$ to be the $k$-module generated by
relative framed links in $(M,\partial M)$ such that $L \cap
\partial M=\partial L=\{x_i, y_i\}$ with the induced framing,
considered up to an ambient isotopy fixing $\partial M$ modulo the
Kauffman skein relations.
\end{Def}

In particular, we will study the relative Kauffman skein module of
the solid torus with two points $A$ and $B$ on the boundary. By a
slight abuse of notation, we denote this module by $K(S^1 \times
D^2, AB)$. Note that $K(S^1 \times D^2, AB)\cong K(S^1 \times D^2,
BA)$ as framed links are not oriented.

We state our main results as the following three theorems.
\begin{Thm}
When $k={\mathbb Q}(\alpha, s)$, $K(S^1 \times D^2, AB)$ has a
countable infinite basis given by the collection of elements of
\{$Q_{\lambda, c}', Q_{\lambda, c}'' $\}, where $\lambda$ varies
over all nonempty Young diagrams and $c$ varies over all extreme
cells of $\lambda$,
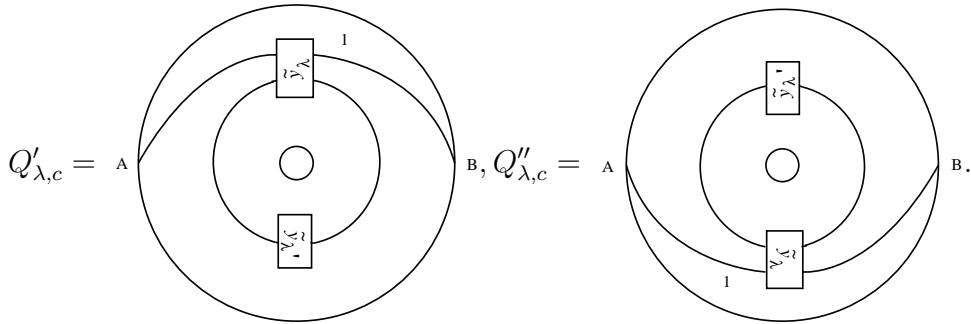
$$
Q_{\lambda, c}'=\raisebox{-20mm}{\
\epsfxsize1.9in\epsffile{relab1.ai}}, Q_{\lambda,
c}''=\raisebox{-20mm}{\ \epsfxsize1.9in\epsffile{relab2.ai}}.
$$
\end{Thm}

\noindent Here the Young diagram $\lambda '$ is obtained from the
Young diagram $\lambda$ by removing the extreme cell $c$. An
{\it extreme cell}  is a cell such that if we remove it, we obtain a
legitimate Young diagram. We draw one single string in the picture
to indicate many parallel strings as to be understood in the
context. A box labelled by $y_{\lambda}$ represents a certain
linear combination of braid diagrams associated to $\lambda$
\cite{AM98} \cite{cB98}. Here the box labelled by
$\widetilde{y_{\lambda}}$ is a homomorphic image of $y_{\lambda}$,
which will be defined in section 2.

\begin{Thm}
When $k={\mathbb Q}(\alpha, s)$, the collection of all the
elements \{$\widehat{\widetilde{y}}_{\lambda}: \lambda$ is any
Young diagram\} forms a basis for the Kauffman skein module of the
solid torus $S^1\times D^2$, which is denoted by $K(S^1\times
D^2)$,
$$\widehat{\widetilde{y}}_\lambda=\raisebox{-24mm}{\
\epsfxsize1.9in
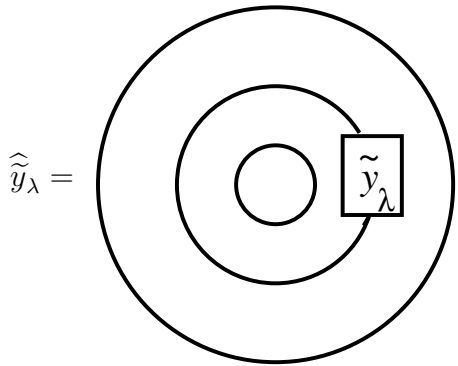
\epsffile{yla1.ai}}.$$
\end{Thm}
\begin{Thm} When $k={\mathbb Q}(\alpha, s)$, $K(S^{1} \times S^2)$
is generated by the empty link, i.e.,
$K(S^{1} \times S^2)=<\phi>.$
\end{Thm}

In section 2, we summarize the work of Beliakava and Blanchet
\cite{BB2000} on the Birman-Murakami-Wenzl category and the
Birman-Murakami-Wenzl algebra $K_n$. In section 3 we compute
$K(S^1\times D^2)$ and prove Theorem 2. In section 4, we study the
relative Kauffman skein module $K(S^1\times D^2, AB)$ and prove
Theorem 1. We use $K(S^1\times D^2, AB)$ to compute
$K(S^{1} \times S^2)$ and prove Theorem 3 in section 5.

\section{The Birman-Murakami-Wenzl category and the
Birman-Murakami-Wenzl algebra $K_n$}

This section is mainly a summary of the related work of  Beliakova
and Blanchet. Details, further references to the origin of some of
these ideas, and related results of others can be found in
\cite{BB2000}. We provide some figures to illustrate the ideas.

\subsection{The Birman-Murakami-Wenzl category and the
Birman-Murakami-Wenzl algebra $K_n$}

{\it The Birman-Murakami-Wenzl category $K$} is defined as: an
object in $K$ is a disc $D^2$ equipped with a finite set of points
and a nonzero vector at each point; if $\beta=(D^2, l_0)$ and
$\gamma=(D^2, l_1)$ are objects, the module $Hom_K(\beta, \gamma)$
is $K(D^2\times [0,1], l_0\times 0 \amalg l_1\times 1)$. For a
Young diagram $\lambda$, we denote by $\square_\lambda$ the object
of the category $K$ formed with one point for each cell of
$\lambda$. We will use the notation $K(\beta, \gamma)$ for
$Hom_K(\beta, \gamma)$. For composition of $f$ and $g$, it's done by stacking $f$
on the top of $g$.
$$K(\beta, \gamma)\times K(\gamma, \delta)\to K(\beta, \delta),$$
$$(f,g)\to fg.$$
Note that Beliakova and Blanchet choose to stack the second one on
the top of the first. Here we follow the convention as in
\cite{ZG99} \cite{AM98}.

As a special case, in the cylinder $D^2\times I$, let $X_n$ be a
set of $n$ distinct framed points on a diameter $D^2\times \{1\}$
and $Y_n$ be a set of $n$ distinct framed points on a diameter
$D^2\times \{0\}$, then the relative Kauffman skein module
$K(D^2\times I, X_n\amalg Y_n)$ is isomorphic to the {\it
Birman-Murakami-Wenzl algebra} $K_n$, which is the quotient of the
braid group algebra $k[B_n]$ by the Kauffman skein relations.

The Birman-Murakami-Wenzl algebra $K_n$ is generated by the
identity ${\mathbf{1}}_n$, positive transpositions $e_1, e_2,
\cdots, e_{n-1}$ and hooks $h_1, h_2, \cdots, h_{n-1}$ as the
following:
$$e_i=\quad \raisebox{-3.5mm}{\epsfxsize.0in\epsffile{ei.ai}}$$
$$h_i=\quad \raisebox{-3.5mm}{\epsfxsize.0in\epsffile{hi.ai}}$$
for $1\leq i\leq n-1$.

Then $K_n$ is the braid group algebra $k[B_n]$ quotient by the
following relations:

\hspace{2cm} $(B_1)\ e_ie_{i+1}e_i=e_{i+1}e_{i}e_{i+1}$,

\hspace{2cm} $(B_2)\ e_ie_j=e_je_{i},\ |i-j|\geq 2$,

\hspace{2cm} $(R_1)\ h_i e_i=\alpha^{-1}h_i$,

\hspace{2cm} $(R_2)\ h_i{e_{i-1}}^{\pm 1}h_i={\alpha}^{\pm 1}h_i$,

\hspace{2cm} $(K)\ e_i-e_i^{-1}=(s-s^{-1})({\mathbf 1}_n-h_i)$.

The quotient of $K_n$ by the ideal $I_n$ generated by $h_{n-1}$ is
isomorphic to the $n$th Hecke algebra $H_n$. Note that
$I_n=\{(a\otimes {\mathbf 1}_1)h_{n-1}(b\otimes {\mathbf 1}_1):\ a,
b\in K_{n-1}\}$. We denote the canonical projection map by
$\pi_n$:
$$\pi_n: K_n\to H_n.$$

\begin{Thm}
There exists a multiplicative homomorphism $s_n:H_n\to K_n$, such
that
$$\pi_n\circ s_n=id_{H_n},$$
$$s_n(x)y=ys_n(x)=0, \ \forall x\in H_n, \forall y \in I_n.$$
\end{Thm}

See details of the proof in \cite{BB2000}.

\begin{Cor}
$K_n\cong H_n\oplus I_n$.
\end{Cor}

\subsection{A basis for the Birman-Murakami-Wenzl algebra
$K_n$}

A sequence $\Lambda=(\Lambda_1, \cdots, \Lambda_n)$ of Young
diagrams will be called an up and down tableau of length $n$ and
shape $\Lambda_n$ if two consecutive Young diagrams $\Lambda_i$
and $\Lambda_{i+1}$ differ by exactly one cell. We observe that in an up
and down tableau $\Lambda=(\Lambda_1, \cdots, \Lambda_n)$ of
length $n$, the size of $\Lambda_n$ is either $n$ or less than $n$
by an even number.

For an up and down tableau $\Lambda$ of length $n$, we denote by
$\Lambda'$ the tableau of length $n-1$ obtained by removing the
last Young diagram in the sequence $\Lambda$. We define
$a_{\Lambda}\in K(n, \square_{\Lambda})$ and $b_{\Lambda}\in K(
\square_{\Lambda}, n)$ by
$$a_1=b_1={\mathbf{1}}_1.$$
If $|\Lambda_n|=|\Lambda_{n-1}|+1$, then
$$a_{\Lambda}=(a_{\Lambda'}\otimes
{\mathbf{1}}_1){\widetilde{y}_{\Lambda_n}},$$
$$b_{\Lambda}={\widetilde{y}_{\Lambda_n}}(b_{\Lambda'}\otimes
{\mathbf{1}}_1);$$ if $|\Lambda_n|=|\Lambda_{n-1}|-1$, then
$$a_{\Lambda}=\frac{<\Lambda_n>}{<\Lambda_{n-1}>}(a_{\Lambda'}\otimes
{\mathbf{1}}_1)({\widetilde{y}_{\Lambda_n}}\otimes \cup),$$
\hspace{3.5cm} $b_{\Lambda}=({\widetilde{y}_{\Lambda_n}}\otimes
\cap)(b_{\Lambda'}\otimes {\mathbf{1}}_1).$

Here $<\lambda>$ is the quantum dimension \cite{W88} associated with
$\lambda$, which is the Kauffman polynomial of
$\widehat{\widetilde{y}}_\lambda$ in $S^3$. $<\lambda>$ is
invertible in ${\mathbb Q}(\alpha, s)$ \cite{BB2000}.

In each case, the figures are drawn below.

(1) If $|\Lambda_n|=|\Lambda_{n-1}|+1$, then
$$a_{\Lambda}=\raisebox{-20mm}{\epsfxsize.0in
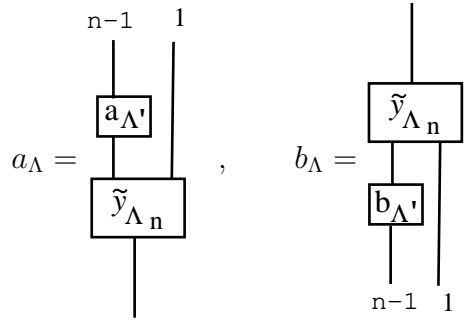
\epsffile{ala.ai}}\quad,\quad  \quad
b_{\Lambda}=\raisebox{-20mm}{\epsfxsize.0in\epsffile{bla.ai}}$$

(2) If $|\Lambda_n|=|\Lambda_{n-1}|-1$, then
$$a_{\Lambda}=\frac{<\Lambda_n>}{<\Lambda_{n-1}>}\quad
\raisebox{-20mm}{\epsfxsize.0in\epsffile{ala1.ai}}\quad, \quad
b_{\Lambda}=\raisebox{-20mm}{\epsfxsize.0in\epsffile{bla1.ai}}$$

\begin{Thm}
The family $a_{\Lambda}b_{\Gamma}$ for all up and down tableaux
$\Lambda, \Gamma$ of length $n$, such that $\Lambda_n=\Gamma_n$
forms a basis for $K_n$.
\end{Thm}

See details of the proof in \cite{BB2000}.

Let $\Lambda=(\Lambda_1, \cdots, \Lambda_n), \
\Gamma=(\Gamma_1,\cdots, \Gamma_n)$ be two up and down tableaux of
length $n$.
If $\Lambda=\Gamma$, i.e. $\Lambda_i=\Gamma_i$ for $1\leq i \leq
n$, then $b_{\Gamma}a_{\Lambda}=\widetilde{y}_{\Lambda_n}$;
otherwise $b_{\Lambda}a_{\Gamma}=0$. This follows from the
corresponding properties in the Hecke category \cite{cB98}. We will use these
properties in the following sections.

\section{The Kauffman skein module of the solid torus $S^1\times
D^2$}

Let $\lambda$ be a Young diagram of size $n$, in the
Hecke category $H_{\square_{\lambda}}$, we have a Young idempotent
$y_{\lambda}$, whose definition and properties are in
\cite[Chapter 3]{ZG99} \cite{cB98}. Let $y_{\lambda}^{*}$ be the
flattened version of $y_{\lambda}$, then $y_{\lambda}^{*}\in H_n$.
Let $\widetilde{y}_{\lambda}=s_n(y_{\lambda}^{*})$ in $K_n$.

The natural wiring of the cylinder $D^2\times I$ into the solid
torus $S^1\times D^2$ induces a homomorphism: $K(D^2\times I, X_n
\amalg Y_n)\to K(S^1\times D^2)$. We denote the image of $K_n$
under the wiring by $\widehat{K_n}$,
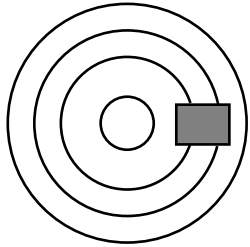
$$\quad \raisebox{-15mm}{\epsfxsize.0in\epsffile{wiring.ai}}.$$

\noindent {\bf Restatement of Theorem 2.} Over ${\mathbb
Q}(\alpha, s)$, the field of rational functions in $\alpha, s$,
the collection of all the elements
\{$\widehat{\widetilde{y}}_{\lambda}: \lambda$ is any Young
diagram\} forms a basis for $K(S^1\times D^2)$.

\begin{proof}
Let $L$ be a framed link in $S^1\times D^2$, up to a scalar multiple, $L$ is the
closure of an $n$-strand braid for some integer $n\geq 0$
\cite{PS}[6.5 Alexander's Braiding Theorem]; the braid modulo the
Kauffman skein relations is an element in $K_n$. Therefore $L\in
\widehat{K_n}$ in $K(S^1\times D^2)$. This shows any element in
$K(S^1\times D^2)$ lies in $\bigcup_{n\geq 0}\widehat{K_n}$, i.e.
$K(S^1\times D^2)\subseteq\bigcup_{n\geq 0}\widehat{K_n}$, hence
$K(S^1\times D^2)=\bigcup_{n\geq 0}\widehat{K_n}$. So we need only
to show that the set \{$\widehat{\widetilde{y}}_{\lambda}:
\lambda$ is a Young diagram\} forms a basis for $\bigcup_{n\geq
0}\widehat{K_n}$.

(1) We first prove that $\widehat{K_n}\subseteq K(S^1\times D^2)$ has a basis given by the
collection \{$\widetilde{y}_{\Lambda_n}: |\Lambda_n|$ is either
$n$ or less than $n$ by an even number \}. It follows that the set 
\{$\widetilde{y}_{\lambda}: 0\leq |\lambda|\leq n$\} forms a basis for $\bigcup_{0\leq i\leq n}\widehat{K_i}$.

By Theorem 5 in the previous section, $K_n$ has a basis given by
the family $a_\Lambda b_\Gamma$, so $\widehat{K_n}$ is generated
by the family $\widehat{a_\Lambda b_\Gamma}$. Since
$\widehat{a_\Lambda b_\Gamma}=\widehat{b_\Gamma
a_\Lambda}=\delta_{\Lambda\Gamma}\widehat{b_\Lambda
a_\Lambda}=\delta_{\Lambda\Gamma}\widehat{\widetilde{y}}_{\Lambda_n}$,
where $\Lambda_n$ is a Young diagram of size either $n$ or less
than $n$ by an even number. so $\widehat{K_n}$ is generated by the
collection \{$\widetilde{y}_{\Lambda_n}: |\Lambda_n|$ is either
$n$ or less than $n$ by an even number \}.

Now we want to prove the above generating set is linearly
independent. We show this by comparing the dimension.

From Corollary 1, we have $K_n\cong H_n\oplus I_n$, so
$\widehat{K_n}\cong \widehat{H_n}+\widehat{I_n}$. Recall that
$I_n=\{(a\otimes {\mathbf 1}_1)h_{n-1}(b\otimes {\mathbf 1}_1):\ a,
b\in K_{n-1}\}$. A typical element in $\widehat{I_n}$ looks like
the following:
$$\quad \raisebox{-3.5mm}{\epsfxsize.0in\epsffile{hatin.ai}}$$
where $a, b \in K_{n-1}$. 
We can see that the above element is in $\widehat{K_{n-2}}$. On the other hand, we
observe that $\widehat{K_{n-2}}\subseteq \widehat{I_n}$, therefore
$\widehat{I_n}=\widehat{K_{n-2}}$. i.e. $\widehat{K_n}\cong 
\widehat{H_n}+\widehat{K_{n-2}}$.

Repeating the process for $\widehat{K_{n-2}}$, we conclude that

\begin{equation}
    \widehat{K_n}\cong
    \begin{cases}
       \widehat{H_1}+\widehat{H_3}+\cdots+\widehat{H_n}, &\textrm {if $n$ is odd;}\notag\\
       <\phi>+\widehat{H_2}+\widehat{H_4}+\cdots+\widehat{H_n},  &\textrm{if $n$ is
       even;}\notag
    \end{cases}
\end{equation}
where $\phi$ denotes the empty link.

Since $\widehat{H_i}\cap \widehat{H_j}=0$ whenever $i\neq j$ in
the Homflypt skein module of the solid torus, (note
$\widehat{H_n}=C_n$ in $S(S^1\times D^2)$ \cite{ZG99}), the above
decomposition is a direct sum.

\begin{equation}
    \widehat{K_n}\cong
    \begin{cases}
       \widehat{H_1}\oplus\widehat{H_3}\oplus\cdots\oplus\widehat{H_n},  &\textrm {if $n$ is odd;}\notag\\
       <\phi>\oplus\widehat{H_2}\oplus\widehat{H_4}\oplus\cdots\oplus\widehat{H_n}, &\textrm{if $n$ is
       even.}\notag
    \end{cases}
\end{equation}
Therefore we have the following equality for the dimensions.

\begin{equation}
    \dim(\widehat{K_n})=
    \begin{cases}
       \dim(\widehat{H_1})+\dim(\widehat{H_3})+\cdots+\dim(\widehat{H_n}),   &\textrm {if $n$ is odd;}\notag\\
       1+\dim(\widehat{H_2})+\dim(\widehat{H_4})+\cdots+\dim(\widehat{H_n}),  &\textrm{if $n$ is even.}
    \end{cases}
\end{equation}

As we know from \cite{AM98} that the dimension of
$\widehat{H_n}=C_n$ in $S(S^1\times D^2)$ is equal to the number
of Young diagrams of size $n$. Therefore the number of generators
in the set \{$\widetilde{y}_{\Lambda_n}: |\Lambda_n|$ is either
$n$ or less than $n$ by an even number \} is equal to the
dimension of $\widehat{K_n}$, therefore it forms a basis for
$\widehat{K_n}$.

(2) The collection of all the
elements \{$\widehat{\widetilde{y}}_{\lambda}: \lambda$ is any
Young diagram\} forms a basis for $K(S^1\times D^2)$.
As we have $K(S^1\times D^2)=\bigcup_{n\geq 0}\widehat{K_n}$.
The result follows by induction on $m$ for $\bigcup_{0\leq n\leq m}\widehat{K_n}$.

\end{proof}

\section{The relative Kauffman skein module $K(S^1\times D^2, AB)$}

\noindent{\bf Restatement of Theorem 1.} When $k={\mathbb
Q}(\alpha, s)$, $K(S^1 \times D^2, AB)$ has a countable infinite
basis given by the collection of elements of \{$Q_{\lambda, c}',
Q_{\lambda, c}'' $\}, where $\lambda$ varies over all nonempty
Young diagrams and $c$ varies over all extreme cells of $\lambda$.
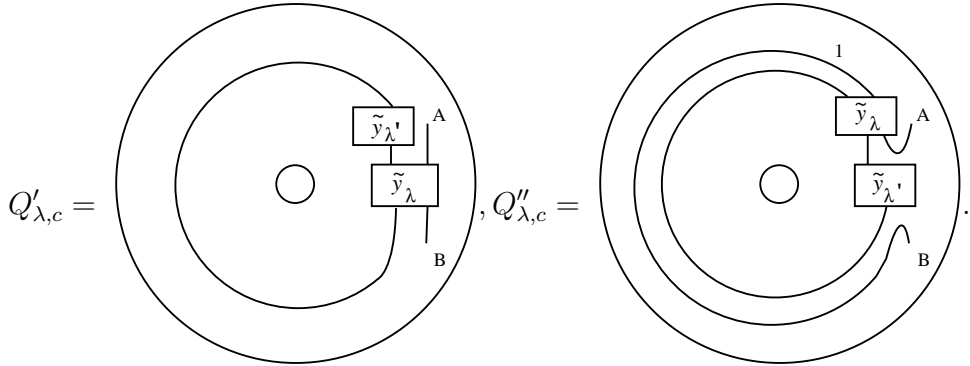
$$
Q_{\lambda, c}'=\raisebox{-20mm}{\
\epsfxsize1.9in\epsffile{baserel1.ai}}, Q_{\lambda,
c}''=\raisebox{-20mm}{\ \epsfxsize1.9in\epsffile{baserel2.ai}}.
$$
Here we change the figures of $Q_{\lambda, c}', Q_{\lambda, c}'' $
through an obvious homeomorphism of $S^1 \times D^2$.

\begin{proof}
(1) $K(S^1 \times D^2, AB)$ is generated by the collection
\{$Q_{\lambda, c}', Q_{\lambda, c}'' $\}.

As we know, $K(S^1 \times D^2, AB)$ is generated by isotopy classes of framed
links in $S^1 \times D^2$ with boundary the two points $A$ and
$B$. Such links consist of a collection of framed closed curves
and a framed arc joining the two points $A$ and $B$. Each framed
link in $S^1\times D^2$ with boundary $A$ and $B$ is isotopic to
one of the two wiring images of an $n$-strand braid for some
$n\geq 1$, which depends on the two ways the last end points of
the top and bottom braid related to $A$ and $B$. The two types of
wirings of the braid group $B_n$ into $K(S^1 \times D^2, AB)$ are
given by the following:

$$
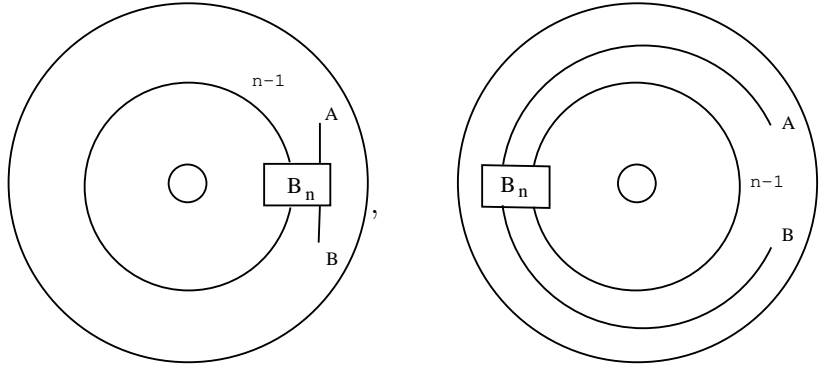
\raisebox{-20mm}{\ \epsfxsize1.9in\epsffile{bnt1.ai}},\quad\quad
\raisebox{-20mm}{\ \epsfxsize1.9in\epsffile{bnt21.ai}}.
$$

By modulo over the Kauffman skein relations, these elements are
linear combinations of elements in the wiring image of
$$
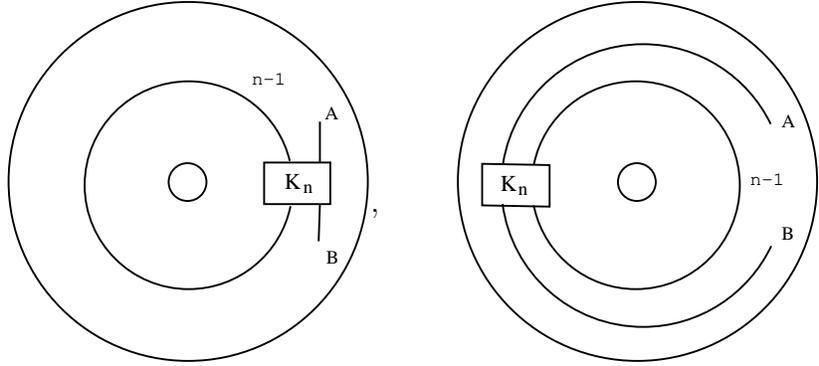
\raisebox{-20mm}{\ \epsfxsize1.9in\epsffile{knt1.ai}},\quad\quad
\raisebox{-20mm}{\ \epsfxsize1.9in\epsffile{knt2.ai}}.
$$
We call the wiring on the left as the first wiring and the wiring
on the right as the second wiring. We denote the images of $K_n$
under the first and second wirings by $\widetilde{K_n}$ and
$\overline{K}_n$, respectively. It is easy to see that
$\widetilde{K_n}=\overline{K}_n$, since $K(S^1\times D^2,
AB)=K(S^1\times D^2, BA)$, which can be seen by
arranging the points $A$ and $B$ symmetrically on the boundary of
$S^1\times D^2$.
Therefore it is sufficient to consider one of the
wiring image. (We will choose the first wiring for convenience.)

In section 2 Theorem 5, $K_n$ has a basis given by the
$a_{\Lambda}b_{\Gamma}$s, therefore the image of any element of
$K_n$ under the first wiring is a linear combination of the wiring
images of $a_{\Lambda}b_{\Gamma}$s given by the following diagrams
according to the two cases that $\Lambda_n$ related to
$\Lambda_{n-1}$ ($|\Lambda_{n-1}|=|\Lambda_{n}|\pm 1$).
$$
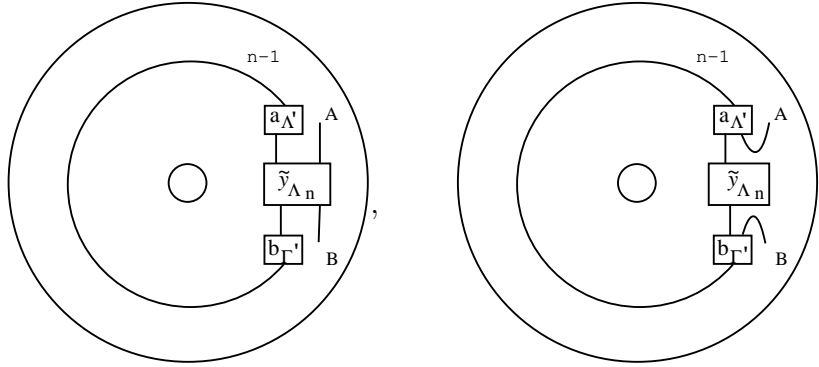
\raisebox{-20mm}{\
\epsfxsize1.9in\epsffile{baserel.ai}},\quad\quad
\raisebox{-20mm}{\ \epsfxsize1.9in\epsffile{baserelr.ai}}.$$ Since
$b_{\Gamma'}a_{\Lambda'}=0$ whenever ${\Gamma'}\neq {\Lambda'}$,
by induction, in the linear combination, the only nonzero elements
are the $b_{\Gamma}a_{\Lambda}$s with ${\Gamma}={\Lambda}$, i.e.
${\Gamma_i}={\Lambda_i}$. Therefore the linear combination will
contain only elements as:
$$
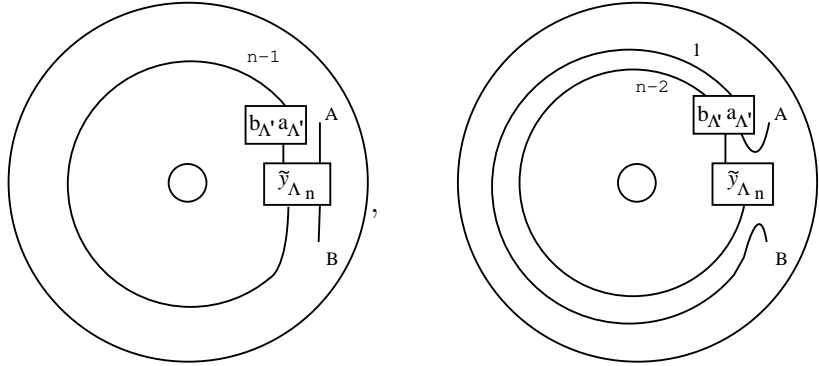
\raisebox{-20mm}{\
\epsfxsize1.9in\epsffile{baserel11.ai}},\quad\quad
\raisebox{-20mm}{\ \epsfxsize1.9in\epsffile{baserel22.ai}}.
$$
Since $b_{\Lambda'}a_{\Lambda'}=\tilde{y}_{\Lambda_{n-1}}$, and
$|\Lambda_{n-1}|=|\Lambda_{n}|\pm 1$, we introduce
$\lambda=\Lambda_{n}$ and $\lambda'=\Lambda_{n-1}$ if
$|\Lambda_{n-1}|=|\Lambda_{n}|- 1$, and introduce
$\lambda=\Lambda_{n-1}$ and $\lambda'=\Lambda_{n}$ if
$|\Lambda_{n-1}|=|\Lambda_{n}|+ 1$; therefore the generators have
the forms
$$
Q_{\lambda, c}'=\quad \raisebox{-20mm}{\
\epsfxsize1.9in
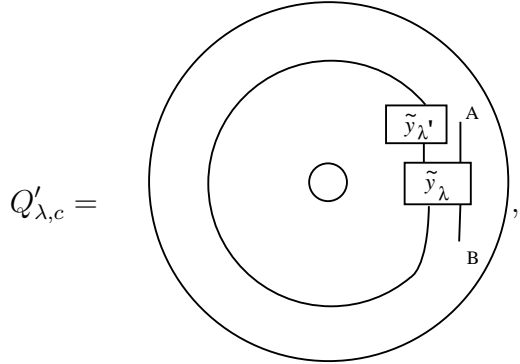
\epsffile{baserel1.ai}},$$
$$Q_{\lambda, c}''=\quad
\raisebox{-20mm}{\ \epsfxsize1.9in
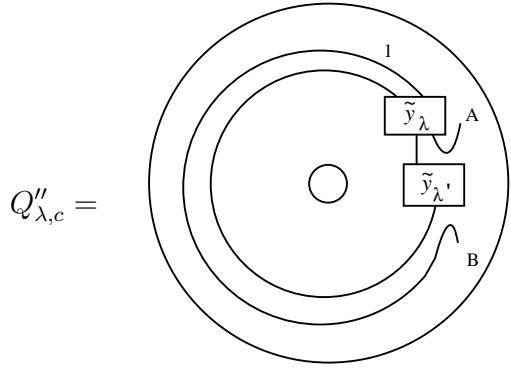
\epsffile{baserel2.ai}},
$$
where $\lambda'$ is a Young diagram obtained from $\lambda$ by
deleting an extreme cell. We conclude that $\widetilde{K_n}$ is generated by the set \{$Q_{\lambda, c}', Q_{\mu, c}''$\} where $|\lambda|\geq 1$, $|\mu|\geq 1$ and $|\lambda|$ is equal to $n$ or less than $n$ by an even number and $|\mu|$ is equal to $n-1$ or less than $n-1$ by an even number.

As $K(S^1\times D^2, AB)=\bigcup_{n\geq 1}\widetilde{K_n}$, the
argument  above shows that the given collection generates
$K(S^1\times D^2, AB)$.

(2) We now prove that the given collection of generators are
linearly  independent. From part (1), we observe that
$\bigcup_{i=1}^{n}\widetilde{K_i}$ is generated by the set
\{$Q_{\lambda, c}', Q_{\mu, c}'' :\ 1\leq |\lambda|\leq n, 1\leq
|\mu|\leq (n-1)$\}. We claim it is sufficient to show that the set
\{$Q_{\lambda, c}', Q_{\mu, c}'' :\ 1\leq |\lambda|\leq n, 1\leq
|\mu|\leq (n-1)$\} forms a basis for
$\bigcup_{i=1}^{n}\widetilde{K_i}$. This is because $K(S^1\times
D^2, AB)=\bigcup_{n\geq 1}\widetilde{K_n}$ and by induction on $m$
for $\bigcup_{1\leq n \leq m}\widetilde{K_n}$, the result follows.

We introduce two types of wirings of $H_n$ into $S(S^1\times
D^2, A, B)$ and denote them by $\widetilde{H_n}$
and $\overline{H}_n$, respectively, where $S(S^1\times D^2, A, B)$
is the relative Homflypt skein modules of $S^1\times D^2$ with $A$
an input point and $B$ an output point on the boundary
\cite{ZG99}[Chapters  2 \& 4].
$$\widetilde{H_n}=\quad
\raisebox{-20mm}{\ \epsfxsize1.9in
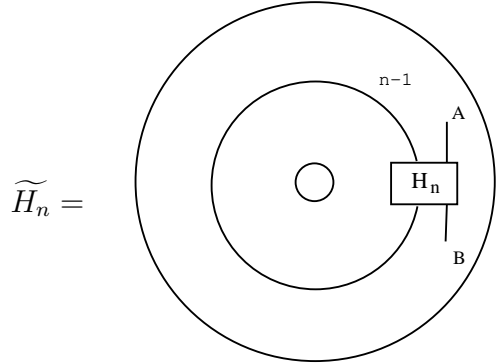
\epsffile{hnt1.ai}}
$$
$$\overline{H}_n=\quad
\raisebox{-20mm}{\ \epsfxsize1.9in
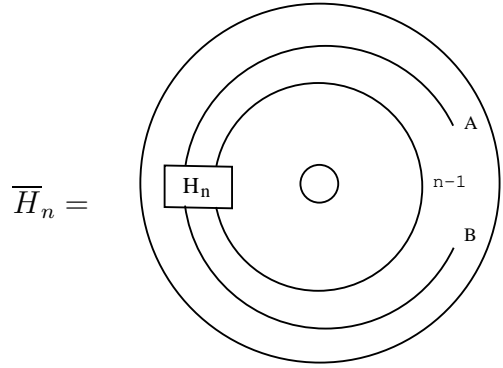
\epsffile{hnt2.ai}}.
$$

From Corollary 1 in section 2, we have $K_n\cong H_n\oplus I_n$, hence
$\widetilde{K_n}\cong \widetilde{H_n}+ \widetilde{I_n}$, where
$\widetilde{I_n}$ is the wiring image of $I_n$ as a submodule of
$K_n$. Now by a similar argument as in the proof of Theorem 1,
$\widetilde{I_n}=\overline{K}_{n-1}$, therefore
$\widetilde{K_n}\cong \widetilde{H_n}+ \overline{K}_{n-1}$,
repeating the process for $\overline{K}_{n-1}$, we have
$$\widetilde{K_n}\cong \widetilde{H_n}+ \overline{H}_{n-1}+\cdots$$
Eventually, we have
$$\bigcup_{i=1}^{n}\widetilde{K_i}\cong \bigcup_{i=1}^{n}\widetilde{H_i}+ \bigcup_{i=1}^{n-1}\overline{H}_i.$$
By the properties of $\widetilde{H_i}$ and $\overline{H}_i$  in
the Homflypt skein module, the above is a direct sum,
$$\bigcup_{i=1}^{n}\widetilde{K_i}\cong \bigcup_{i=1}^{n}\widetilde{H_i}\oplus \bigcup_{i=1}^{n-1}\overline{H}_i.$$

Therefore, we have the following equality of the dimensions,

$$\dim(\bigcup_{i=1}^{n}\widetilde{H_i})+\dim(\bigcup_{i=1}^{n-1}\overline{H}_i)= \dim(\bigcup_{i=1}^{n}\widetilde{K_i})$$

Note that $\widetilde{H_i}$ is the subspace $C_i'$ \cite[Chapter
4]{ZG99} in the relative Homflypt skein module $S(S^1\times D^2,
A, B)$. A basis of $C_i'$ was given in \cite[Chapter 4]{ZG99} as
\{$Q_{H_{\lambda, c}'}:\ |\lambda|=i $\},
$$
Q_{H_{\lambda, c}'}=\raisebox{-33mm}{\
\epsfxsize0in\epsffile{hnrel1.ai}}.
$$

Similarly, a basis of $\overline{H}_{i}=C_{-i}'$ was given  in
\cite[Chapter 4]{ZG99} as \{$Q_{H_{\lambda, c}''}:\ |\lambda|=i
$\}.
$$Q_{_H{\lambda, c}''}=\raisebox{-33mm}{\
\epsfxsize0in\epsffile{hlacc.ai}}.$$

Therefore the set \{$Q_{H_{\lambda, c}'}, Q_{H_{\mu, c}''} :\ 1\leq
|\lambda|\leq n, 1\leq |\mu|\leq (n-1)$\}  forms a basis for
$\bigcup_{i=1}^{n}\widetilde{H_i}\oplus
\bigcup_{i=1}^{n-1}\overline{H}_i$.
Hence the cardinality of the above set is
$\dim(\bigcup_{i=1}^{n}\widetilde{H_i})+\dim(\bigcup_{i=1}^{n-1}\overline{H}_i)$.

On the other hand, since $\bigcup_{i=1}^{n}\widetilde{K_i}$ is
generated by the set \{$Q_{\lambda, c}', Q_{\mu, c}''\  : 1\leq
|\lambda|\leq n, 1\leq |\mu|\leq (n-1)$\}. As
this set has the same cardinality as the set \{$Q_{H_{\lambda,
c}'}, Q_{H_{\mu, c}''} :\ 1\leq |\lambda|\leq n, 1\leq |\mu|\leq
(n-1)$\}, with
$$\dim(\bigcup_{i=1}^{n}\widetilde{H_i})+\dim(\bigcup_{i=1}^{n-1}\overline{H}_i)= \dim(\bigcup_{i=1}^{n}\widetilde{K_i})$$ from above, we conclude that the set  \{$Q_{\lambda, c}', Q_{\mu, c}'' :\ 1\leq
|\lambda|\leq n, 1\leq |\mu|\leq (n-1)$\} must be a basis for
$\bigcup_{i=1}^{n}\widetilde{K_i}$.

\end{proof}

\section{The Kauffman skein module $K(S^1\times S^2)$}

The space $S^1 \times S^2$ can be obtained from the solid torus
$S^1 \times D^2$ by first attaching a $2$-handle along the
meridian $\gamma$ and then attaching a $3$-handle. As it is
well-known, adding a $3$-handle induces isomorphism between skein
modules; while adding a $2$-handle to $S^1\times D^2$ adds
relations to the module ${K(S^1\times D^2)}$. The natural
inclusion $i$: ${S^1\times D^2} \to {S^1\times S^2}$ induces an
epimorphism $i_{*}$: ${K(S^1\times D^2)} \to {K(S^1\times S^2)}$.
Following Masbaum's work in the case of the Kauffman bracket skein
module \cite{gM96}, we will use the following method to
parametrize the relations arising from sliding over the
$2$-handle. Pick two points $A, \ B$ on $\gamma$, which decompose
$\gamma$ into two intervals $\gamma '$ and $\gamma ''$.

By Corollary 2 in  \cite[Chapter 2]{ZG99} , we have

$$K(S^1\times S^2)\cong  K(S^1\times D^2)/ R.$$

\noindent Here $R=\{\Phi '( z ) - \Phi ''( z ) \mid z \in K(S^1
\times D^2, AB)\}$, $z$ is any element in the relative skein
module $K(S^1 \times D^2, AB)$, and $\Phi '( z )$ and $\Phi ''( z
)$ are given by capping off $z$ with $\gamma '$ and $\gamma ''$,
respectively, and pushing the resulting links back into $S^1
\times D^2$.

In the previous section we give  a basis for the relative skein
module $K(S^1 \times D^2, AB)$. We now compute
$K(S^1\times S^2)$.

\begin{Thm}
The following is a complete set of relations:
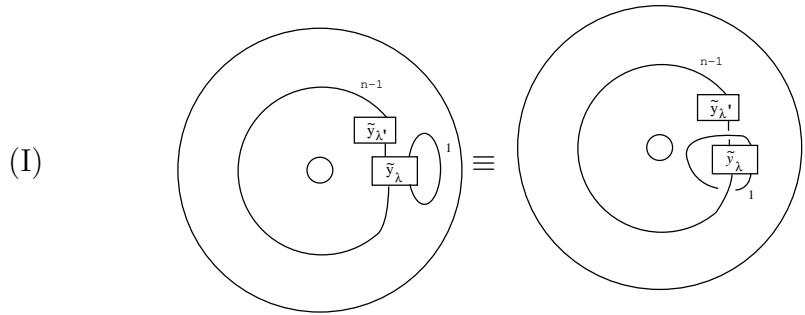
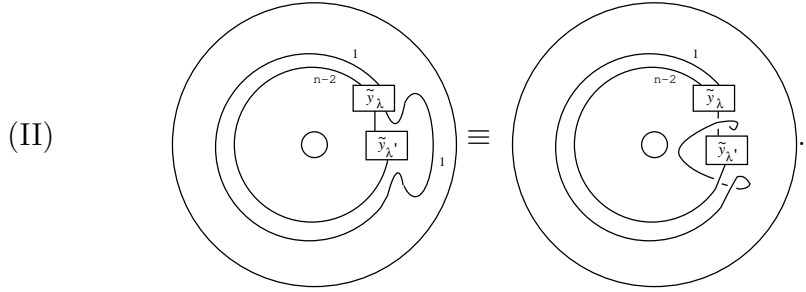
\begin{equation}
\raisebox{-19mm}{\ \epsfxsize1.5in\epsffile{baser1.ai}} \equiv
\raisebox{-16mm}{\ \epsfxsize1.5in\epsffile{baser12.ai}} \tag{I}
\end{equation}
\begin{equation}
\raisebox{-19mm}{\ \epsfxsize1.5in\epsffile{baser2.ai}} \equiv
\raisebox{-19mm}{\ \epsfxsize1.5in\epsffile{baser22.ai}}. \tag{II}
\end{equation}
\end{Thm}

\begin{proof}
This follows from $\Phi '( z ) \equiv \Phi ''( z )$ by taking $z$
to be every basis element in $K(S^1 \times D^2, AB)$, which gives
a generating set for the submodule $R\subset K(S^1\times D^2)$.
\end{proof}

\noindent{\bf Restatement of Theorem 3.} Over $k={\mathbb
Q}(\alpha, s)$, $K(S^{1} \times S^2)=<\phi>.$

\begin{proof}
We show that every basis element $\widehat{\widetilde{y}}_\lambda$
of $K(S^{1} \times D^2)$ with $\lambda$ being nonempty is $ 0$ in
$K(S^{1} \times S^2)$. We will use the complete set of relations
given in the previous theorem to compute $R$.

(1) We simplify equation (I) in Theorem 6. The left hand side is equal to
a scalar multiple of $\widehat{\widetilde{y}}_{\lambda'}$, by
embedding it into $S^3$, the scalar is
$\frac{<\lambda'>}{<\lambda>}$. Here ${<\lambda>}$ is the Kauffman
polynomial of $\widehat{\widetilde{y}}_\lambda$ in $S^3$. We can
simplify the right hand side of (I) by using the following skein
relation \cite[Prop. 6.1]{BB2000}:
$$\raisebox{-27mm}{\ \epsfxsize0in\epsffile{alasb.ai}}
=s^{2cn(c)}\quad \raisebox{-12mm}{\
\epsfxsize0in\epsffile{alasa.ai}},$$
where $cn(c)$ is the content
of the extreme cell $c$ of $\lambda$ to be removed to obtain
$\lambda'$.

The right hand side of (I) is also a scalar multiple of
$\widehat{\widetilde{y}}_{\lambda'}$, by embedding it into $S^3$ and using the skein relation
above, the scalar is equal to
$s^{2cn(c)}(\frac{<\lambda'>}{<\lambda>})$. Therefore, equation
(I) gives the equivalence relation:
$\frac{<\lambda'>}{<\lambda>}\widehat{\widetilde{y}}_{\lambda'}\equiv
s^{2cn(c)}(\frac{<\lambda'>}{<\lambda>})\widehat{\widetilde{y}}_{\lambda'}
$, i.e. $\frac{<\lambda'>}{<\lambda>}(1-s^{2cn(c)})
\widehat{\widetilde{y}}_{\lambda'}\equiv 0$ in $K(S^1\times S^2)$.
Since both $\frac{<\lambda'>}{<\lambda>}$ and $(1-s^{2cn(c)})$ are
invertible, we get $\widehat{\widetilde{y}}_{\lambda'}\equiv 0$ in
$K(S^1\times S^2)$ when $\lambda'$ is not empty. 

(2) Now we simplify equation (II) in Theorem 6. The left hand side is equal to
$\widehat{\widetilde{y}}_\lambda$ by the absorbing property. The
right hand side is equal to
$\alpha^{-2}s^{2cn(c)}\widehat{\widetilde{y}}_\lambda$ by the
skein relation above. So equation (II) gives $
(1-\alpha^{-2}s^{2cn(c)})\widehat{\widetilde{y}}_\lambda\equiv 0$,
which also implies  that
$\widehat{\widetilde{y}}_\lambda\equiv 0$ over $\mathbf{Q}(\alpha, s)$, so equation (II)  gives no
new information.

Therefore all $\widehat{\widetilde{y}}_\lambda\equiv 0$ in
$K(S^{1} \times S^2)$ when $\lambda$ is nonempty.

(3) As no relation involves the empty link $\phi$, it
survives. Therefore,

$$K(S^{1} \times S^2)=<\phi>.$$

\end{proof}

\end{document}